\documentclass[12pt]{article}
\usepackage{amsmath}
\usepackage{amsfonts}
\usepackage{amsthm}
\usepackage{amssymb}
\usepackage{color}
\usepackage{tkz-graph}

\usepackage{graphicx}
\usepackage{listings}
\usepackage{tkz-berge}

\newtheorem{theorem}{Theorem}[section]
\newtheorem{lem}[theorem]{Lemma}

\newtheorem{problem}[theorem]{Problem}
\newtheorem{cor}[theorem]{Corollary}

\theoremstyle{definition}
\newtheorem{definition}[theorem]{Definition}
\theoremstyle{remark}
\newtheorem{rem}[theorem]{Remark}

\numberwithin{equation}{section}
\numberwithin{theorem}{section}
\numberwithin{figure}{section}

\def\f2{\mathbb{F}_2}

\newcommand{\diam}{{\rm diam}\hskip0.02cm}

\DeclareMathOperator{\Lam}{\mathsf{La}}

\newcommand{\bbN}{\mathbb{N}}
\newcommand{\bn}{\mathbb{N}}

\newcommand{\al}{\alpha}
\newcommand{\be}{\beta}
\newcommand{\g}{\gamma}

\newcommand{\e}{\varepsilon}

\newcommand{\vt}{\theta}
\newcommand{\la}{\lambda}

\newcommand{\vf}{\varphi}

\newcommand{\bbZ}{\mathbb{Z}}

\newcommand\remove[1]{}

\newcommand{\lb}{\label}

\newcommand{\lra}{\longrightarrow}

\newcommand{\wtw}{if and only if}

\newcommand{\Buo}{Without loss of generality }

\newcommand{\DEF}{\buildrel {\mbox{\tiny def}}\over =}

\newcommand{\tsp}{\operatorname{tsp}}

\newcommand{\mD}{\mathsf{D}}
\newcommand{\WB}{\mathsf{WB}}
\newcommand{\LWB}{\mathsf{LWB}}
\newcommand{\symdif}{\bigtriangleup} 

\begin{document}

\title{Bi-Lipschitz embedding properties of lamplighter graphs on weighted and unweighted  trees}

\author{Charlotte Melby\thanks{Partially funded by the Miami University Undergraduate Summer Scholars Program  2024} \ and Beata~Randrianantoanina}

\date{}
\maketitle

\begin{abstract}
In 2021 Baudier, Motakis, Schlumprecht, and Zs\'ak  proved that if a sequence of graphs
 $(G_k)_{k\in{\mathbb{N}}}$ contains the sequence of complete graphs with uniformly bounded distortion, then the sequence of lamplighter graphs on $G_k$'s  contains Hamming cubes with 
uniformly bounded distortion and  asked whether the converse holds. 

They suggested 
that a sequence of trees with edges replaced by paths of ``moderately growing'' lengths may be a counterexample. 
We prove that indeed this is the case, and that   a sequence of
``moderately'' weighted trees is another counterexample.

Further, we prove that  diamond graphs do not embed with uniformly bounded distortion into  lamplighter graphs on 
 trees with edges replaced by paths with sufficiently fast growing  lengths.
\end{abstract}

\smallskip



\section{Introduction}

Given metric spaces $(M,d_M)$,  $(N,d_N)$, and a constant $C\ge 1$, a map $f:M\lra N$ is called a {\it bi-Lipschitz embedding with distortion at most $C$} if there exists a scaling factor   $\la>0$ so that for all
 $u,v\in M$
\begin{equation}
\la d_M(u,v)\le d_N(f(u),f(v))\le C\la d_M(u,v).
\end{equation}
We denote
$$c_N(M)\DEF \inf\{ dist(f):  f:M\to N \text{ is a bi-Lipschitz embedding} \}.$$

If $(M_k)_{k\in\bn}$ and $(N_k)_{k\in\bn}$ are sequences of metric
spaces, we say that \emph{$(M_k)_{k\in\bn}$ equi-bi-Lipschitzly embeds
  into $(N_k)_{k\in\bn}$} if there exists a constant 
$C>0$ such that for all $k\in\bn$ there exists $n\in\bn$ such that
$c_{N_n}(M_k)\le C$.

The theory of bi-Lipschitz and metric embeddings has many interactions with  geometry, analysis, probability, combinatorics, group theory, topology, and theoretical computer science. It has been studied intensively over the last four decades, see \cite{Nao12,Ost13},  and continues to be  a very active and growing area of research, see e.g. the recent workshop on {\it Metric embeddings} at  the American Institute of Mathematics \cite{AIM2025}. 
In this theory 
one  aims to study a structure of an object by analyzing its embedding properties with respect to simpler metric or geometric objects, such as trees, which are  one of the simplest, yet very rich and important objects, both in the geometric group theory and in computer science (cf. \cite{LNP09} and its references). Here, by a tree we mean a (weighted or unweighted) connected graph $T$ that does not contain any cycles, i.e., such that for any $x,y\in T$ there exists exactly one shortest path from $x$ to $y$.
Trees    play a fundamental role in hyperbolic geometry and also  in computer science, see \cite{LNP09} and its references.
 The bi-Lipschitz structure of trees has been studied extensively.
In particular, in a seminal 
paper \cite{Bou86} Bourgain proved
 that  the sequence  of finite complete unweighted binary trees $B_k$ of depth $k$ equi-bi-Lipschitzly  embeds 
into a 
 Banach space $X$ \wtw\ the space $X$ is not 
 superreflexive (i.e. does not admit an equivalent uniformly convex norm), including $X=\ell_1$. 
This 
characterization of a linear property of Banach spaces in terms of their metric structure was one of the influential results that gave rise to  the Ribe program and  started  
a considerable amount of work  in related directions, cf. \cite{Nao12,Ost13}. Notably, Johnson and Schechtman \cite{JS09} 
proved  that the sequence of diamond graphs  $(\mD_k)_{k\in\bbN}$
 (see Definition~\ref{defdiamond} in Section~\ref{sec-pre}) equi-bi-Lipschitzly embeds into
a Banach space $X$ if and only if    $X$ is not   superreflexive.

 In the present paper we are interested in  lamplighter graphs, whose  systematic study was 
initiated by Baudier, Motakis, Schlumprecht, and Zs\'ak in  \cite{BMSZ21}.
Given a graph $G=(V(G),E(G))$, the {\it lamplighter graph} $\Lam(G)$ is the graph whose vertices are of the form 
$(A,x)$, where $x$ is a vertex of $G$ and $A$ is a finite subset of $V(G)$. Two vertices $(A,x), (B,y)$ of 
$\Lam(G)$ are connected by an edge \wtw\ $(i)$ $A=B$ and $(x,y)$ is an edge in $G$, or $(ii)$ $x=y$ and $A\symdif B=\{x\}$, where $A\symdif B$ denotes the symmetric difference of the sets $A, B$, i.e. 
$A\symdif B=(A\setminus B)\cup(B\setminus A)$, (see \cite{BMSZ21}). 
Intuitively, a vertex $(A,x)\in \Lam(G)$ ``consists'' of the set $A\subseteq V(G)$ that represents the set of lamps that are ``on'' and $x\in V(G)$ represents the current position of the lamplighter. The lamplighter can make one of two types of moves: either move to vertex connected by an edge to the current position, or switch the light ``on''/``off'' at the current position, cf. \cite{T17}. By \cite[Proposition~2.1]{BMSZ21}, the metric on $\Lam(G)$ is given by
\begin{equation}\lb{tsp}
d_{\Lam(G)}((A,x),(B,y))=\tsp_G(x,A\symdif B,y) +|A\symdif B|,
\end{equation}
where $|S|$ denotes the cardinality of the set $S$, and $\tsp_G(x,A\symdif B,y)$ denotes the solution of the traveling salesman problem, that is, the 
length of the shortest possible walk on $G$ that starts at $x$, ends at $y$, and passes through all vertices $v\in A\symdif B$.

This definition of lamplighter graphs  coincides  with the standard definition of the wreath product 
$ \bbZ_2 \wr G$ in the case when $G$ is  (a Cayley graph of) a group and it is  its very natural generalization to the setting of arbitrary graphs (not necessarily Cayley graphs).  Since their introduction, the lamplighter graphs have already been used to prove results about the lamplighter groups themselves, see \cite{G22,GT24}.

Lamplighter groups 
are a very interesting class of groups that has been a rich source of important examples and counterexamples in geometric group theory, see \cite{B17,G22,GT24} and their references. In particular, the study of embedding properties of wreath products into the Hilbert space $L_2$ and into $L_1$ has received a lot of attention in the last 20 years -- we refer interested readers to an excellent overview of the literature in \cite{NP11}, cf. also \cite{NP08,G22,BMSZ21,BMSZ22} and their references. 
Here we mention only a few results most closely related to the topic of the present paper.

Lyons, Pemantle, and Peres \cite{LPP96} (cf.  \cite{OR19} and \cite{BMSZ21} for alternate short proofs) proved that the sequence of complete binary trees $(B_k)_k$ 
equi-bi-Lipschitzly embeds into the lamplighter group   $ \bbZ_2 \wr \bbZ$  (and therefore, by the characterization of Bourgain mentioned above, $ \bbZ_2 \wr \bbZ$ does not bi-Lipschitzly embed into any superreflexive Banach space).
 Naor and Peres in \cite{NP08} gave two proofs that  the sequence of  wreath products  $(\bbZ_2 \wr \bbZ_k)_{k\in\bbN}$ equi-bi-Lipschitzly embeds into $\ell_1$, and so does 
 $ \bbZ_2 \wr \bbZ$ (here $\bbZ_k$ denotes the cyclic group of order $k$).  In  \cite{OR19} it was shown that  the sequence   $(\bbZ_2 \wr \bbZ_k)_{k\in\bbN}$  and $ \bbZ_2 \wr \bbZ$ equi-bi-Lipschitzly embed into every nonsuperreflexive Banach space, thus identifying  another metric characterization of superreflexivity.
Cornulier, Stalder, and Valette \cite{CSV12}  showed that for every finitely generated free group $F$, the wreath product $\bbZ_2 \wr F$, and therefore the lamplighter group $\Lam(T)$ on any tree $T$, bi-Lipschitzly embeds into $\ell_1$,  see also \cite{BMSZ21} for a  
purely metric, 
constructive proof.

The starting point for the present work is the paper \cite{BMSZ21} of Baudier, Motakis, Schlumprecht, and Zs\'ak. One of the main results of  \cite{BMSZ21} is
\begin{theorem}{\rm 
 \cite[Theorem~1.3]{BMSZ21}
}\lb{thmC}
The sequences 
$(\Lam(B_k))_{k\in \bbN}$, $(H_k)_{k\in \bbN}$, and $(K_k)_{k\in \bbN}$ pairwise equi-bi-Lipschitzly embed into each other,
 where,    
for   $k\in\bbN$, $K_k$ is a complete graph on $k$ vertices and
$H_k$ is the Hamming cube on the set $\{1,\dots,k\}$ (for full definition see Section~\ref{sec-pre} below).
\end{theorem}

As one of the tools for proving Theorem~\ref{thmC}, the authors of  \cite{BMSZ21}   showed 

\begin{theorem}{\rm \cite[Lemmas~6.2 and~5.1]{BMSZ21}} \lb{complete}
If a sequence of complete graphs $(K_k)_k$ equi-bi-Lipschitzly embeds
  into a sequence of graphs $(G_k)_k$, then the sequence of Hamming cubes $(H_k)_k$ equi-bi-Lipschitzly embeds
  into $(\Lam(G_k))_k$.
\end{theorem}

The authors of  \cite{BMSZ21}    asked whether the converse of Theorem~\ref{complete} is  true.

\begin{problem}{\rm \cite[Problem~7.1]{BMSZ21}}
  \label{problem:7.1}
  Given a sequence $(G_k)_{k\in\bbN}$ of graphs such that the Hamming cubes
  $(H_k)_{k\in\bbN}$ equi-bi-Lipschitzly embed into
  $(\Lam(G_k))_{k\in\bbN}$, does it follow that $(K_k)_{k\in\bbN}$
  equi-bi-Lipschitzly embed into $(G_k)_{k\in\bbN}$?
\end{problem}

They
suggested that a counterexample might be a sequence of 
trees $W_{2,n}$, obtained by replacing each edge of level $k$ of a binary tree $B_n$  by a path of length $2^{n-k}$, as  described in 
\cite[ Example in Section~7]{BMSZ21},  see Section~\ref{ssec-Wlan} below for a precise definition.
We prove that, indeed, this construction 
leads to a counterexample for  Problem~\ref{problem:7.1} and that another counterexample is 
 a sequence of
``moderately'' weighted complete binary trees $(\WB_{2,n})_n$, see  Section~\ref{ssec-Wlan} below. Namely, we prove   
\begin{theorem}\lb{main1}
 For any $\vt\in\bbN$ with $\vt\ge 2$,  the sequence of complete graphs $(K_k)_{k\in\bbN}$ does not   equi-bi-Lipschitzly embed  into the 
sequence of  trees  $(W_{\vt,n})_{n\in\bbN}$ or into the sequence of  weighted trees  
$(\WB_{\vt,n})_{n\in\bbN}$.
\end{theorem}
\begin{theorem}\lb{main2}
 For every $\e>0$, the sequence of Hamming cubes
  $(H_k)_{k\in\bbN}$ bi-Lipschitzly embeds  into $(\Lam(W_{2,n}))_{n\in\bbN}$ and into
 $(\Lam(\WB_{2,n}))_{n\in\bbN}$ with distortion at most $(1+\e)$.
\end{theorem}

As pointed out in \cite{BMSZ21}, since 
there exist nonsuperreflexive Banach spaces that do not  equi-bi-Lipschitzly  contain 
the sequence of Hamming cubes, see  \cite{BMW86,James74,James78,Pisierbook}, 
the above counterexamples also provide  a negative answer to  \cite[Problem~7.3]{BMSZ21}. 
\begin{cor} The sequences of trees  $(W_{2,n})_{n}$  and $(\WB_{2,n})_{n}$
do not
   contain $(K_k)_{k}$ equi-bi-Lipschitzly,
but  both sequences of lamplighter graphs on them 
 $(\Lam(W_{2,k}))_{k\in\bbN}$  and  $(\Lam(\WB_{2,n}))_{n}$ equi-bi-Lipschitzly contain the sequence of Hamming cubes, and thus
do not  equi-bi-Lipschitzly embed into every  nonsuperreflexive Banach space.
\end{cor}

We note that  both sequences   $(W_{2,n})_{n}$  and $(\WB_{2,n})_{n}$  embed equi-bi-Lipschitzly into the Hilbert space $L_2$.
This follows from \cite[Theorem~1.1]{LNP09} which says 
that a metric tree admits a bi-Lipschitz embedding into Hilbert space if and only if
it does not equi-bi-Lipschitzly contain the sequence of complete binary trees  $(B_n)_{n}$. Indeed, 
since the sequence of complete graphs $(K_k)_k$ equi-bi-Lipschitzly embeds into the sequence 
$(B_n)_{n}$ (\cite{Mat99}), 
Theorem~\ref{main1} implies that  
$(B_n)_{n}$ does not equi-bi-Lipschitzly embed  into either of
 $(W_{2,n})_{n}$ or $(\WB_{2,n})_{n}$ 
(alternatively, one can  note that the  trees $(W_{2,n})_{n}$ and 
$(\WB_{2,n})_{n}$ are doubling and use an earlier result \cite[Theorem~2.4]{GKL03}).

Thus sequences   $(W_{2,k})_{k}$  and $(\WB_{2,n})_{n}$ are new examples of graphs, other than $(K_n)_{n}$,
that equi-bi-Lipschitzly embed into Hilbert space and such that the  lamplighter graphs on them do not  equi-bi-Lipschitzly embed into every  nonsuperreflexive Banach space.

We  note that, as is easy to see,   the sequence   $(W_{2,k})_{k}$  embeds equi-bi-Lipschitzly into the sequence of diamonds and also into the sequence of Laakso graphs. Thus,  by \cite[Lemma~5.1]{BMSZ21}, 
Theorem~\ref{main2} implies the following corollary:
\begin{cor} 
The  sequence of Hamming cubes embeds equi-bi-Lipschitzly into the sequence of lamplighter graphs on the diamond graphs and also  into the sequence of lamplighter graphs on the Laakso graphs.
\end{cor}

Next we prove that  Theorem~\ref{main2} fails for $(\Lam(W_{\vt,n}))_{n}$ when $\vt\ge 3$ in the following much stronger sense (recall that the sequence of diamond graphs  $(\mD_k)_{k}$ equi-bi-Lipschitzly embeds into the sequence of Hamming cubes
  $(H_k)_{k}$, \cite{GNRS04}, cf. \cite{JS09}).

\begin{theorem}\lb{nodiamonds}
For all $\vt\ge 3$, the sequence of diamonds $(\mD_k)_{k\in\bbN}$ does not equi-bi-Lipschitzly embed into the  sequence of  lamplighter graphs  $(\Lam(W_{\vt,n}))_{n\in\bbN}$.
\end{theorem}

Since $(W_{\vt,n})_{n}$ are trees, it follows from  \cite{CSV12} (cf. \cite[Theorem~1.1]{BMSZ21}) that the sequence $(\Lam(W_{\vt,n}))_{n}$ embeds equi-bi-Lipschitzly into $\ell_1$. 
However, we do not know whether 
$(\Lam(W_{\vt,n}))_{n}$ when $\vt\ge 3$, embed  equi-bi-Lipschitzly 
into every nonsuperreflexive Banach space. If yes,  the sequences $(\Lam(W_{\vt,n}))_{n}$ for 
$\vt\ge 3$, would
 be a new
metric characterization of superreflexivity. If no, this would be an even more interesting example.

We do not know whether Theorem~\ref{nodiamonds} 
is true for the sequence of lamplighter graphs
on the weigthed complete binary trees $(\Lam(\WB_{\vt,n}))_{n}$, for all or some $\vt\ge 3$.
\smallskip

\noindent
{\bf Acknowledgements.} The second named author is greatful to Chris Gartland for pointing out, during a recent workshop at AIM,  the result of Alexandros Eskenazis that finite products of trees do not contain diamonds uniformly 
(see Corollary~\ref{AE} below), and to the organizers for their invitation to participate in the Workshop on ``Metric Embeddings" at the American Institute of Mathematics. 

   

\section{Preliminaries}\lb{sec-pre}

We use standard notation as may be found, e.g. in \cite{LT77,Ost13}.

Recall that a graph $G=(V,E)$ is a pair of the set of vertices $V=V(G)$ and the set of edges $E=E(G)\subseteq V\times V$. Frequently, with a slight abuse of notation, we will write  that a vertex $v\in G$ instead of writing $v\in V(G)$. 
A {\it walk} in $G$  is a sequence $(x_0,x_1,\dots,x_m)$ of vertices of $G$ such that for every $i$ with $1\le i\le m$,
$(x_{i-1},x_i)$ is an edge of $G$. We will only consider connected graphs, that is, graphs such that for every $x,y\in G$ there exists at least one walk in $G$ from $x$ to $y$.
The unweighted connected graph $G$ is endowed with the shortest path metric $d_G$, where for all $x,y\in G$
\[d_G(x,y)=\min\{m : \text{there exists a walk of length $m$ from $x$ to $y$}\}.\]

A walk  from $x$ to $y$ is called a {\it (geodesic) path}, denoted  $p(x,y)$, if its  length is equal to the graph distance between $x$ and $y$.

Our main interest here is in graphs known as trees, Hamming cubes, and diamonds. 

Recall that a graph $G$ is called a {\it tree}, if $G$ is connected and does not contain any cycles, i.e., if for any $x,y\in G$ there exists exactly one path $p(x,y)$.

For an arbitrary set $I$, the {\it Hamming cube $H_I$} is a graph with the vertex set $\{0,1\}^{(I)}$ consisting of all functions 
$x:I\to\{0,1\}$ with finite support. 
Two vertices 
$x,y\in H_I$ are joined by an edge if and only if they differ in exactly one coordinate, i.e., there exists a unique $i\in I$ with $x_i\neq y_i$. The graph distance on $H_I$ is the $\ell_1$ metric given by
\[d_H(x,y)=\sum_{i\in I}|x_i-y_i|\].

\begin{definition}{\rm(cf. \cite{Ost13})} \lb{defdiamond}
The {\it diamond graphs $\mD_k$}  ($k\in \bbN$) are constructed inductively as follows: $\mD_0$ is a single edge. 
If the graph $\mD_{k-1}$ is already constructed, the diamond $\mD_k$ is obtained by replacing each edge 
$(u,v)$ of $\mD_{k-1}$ by a a quadrilateral $u,a,v,b$, with edges $(u,a)$, $(a,v)$, $(v,b)$,  $(b,u)$. Equivalently, for any $m<k$, one may think of a diamond $\mD_k$ as obtained by replacing each edge of a diamond $\mD_m$ by a diamond 
$\mD_{k-m}$.  Each (isometric)  copy of $\mD_{k-m}$ inside $\mD_k$ that replaced an edge of $\mD_m$, is called
a {\it subdiamond} of $\mD_k$. 

We denote the graph metric on $\mD_k$ by $d_\mD$.
\end{definition}

The following lemma 
follows directly from the definition of the diamond graphs. 

\begin{lem} \lb{subdiamond}
Let $k,j\in\bbN$ with  $k>j+3$, and $u,v\in \mD_k$. If $d_\mD(u,v)\ge 2^{k-j}$, then there exists a subdiamond  $\bar{\mD}=\mD_{k-j-3}\subset \mD_k$ 
that lies between $u$ and $v$, that is, such that for all $w\in\bar{\mD}$,
$d_\mD(u,w)+d_\mD(w,v)=d_\mD(u,v)$.
\end{lem}

\section{Definitions of the 
trees $W_{\vt,n}$ and
$\WB_{\vt,n}$}\lb{ssec-Wlan}

Given   $n\in\bbN$ and $\vt\in\bbN$, $\vt\ge 2$, to define the 
tree $W_{\vt,n}$, we start with a rooted binary tree $B_n$. We say that a vertex of $B_n$ is on the {\it level $k$}, $0\le k \le n$, if its distance from the root is equal to $k$, and that 
an edge of $B_n$ is on the level $k$, $1\le k\le n$, if its lower (that is,  further from the root) vertex  is on the level $k$.

\begin{figure}
    \centering
    \includegraphics[width=0.5\linewidth]{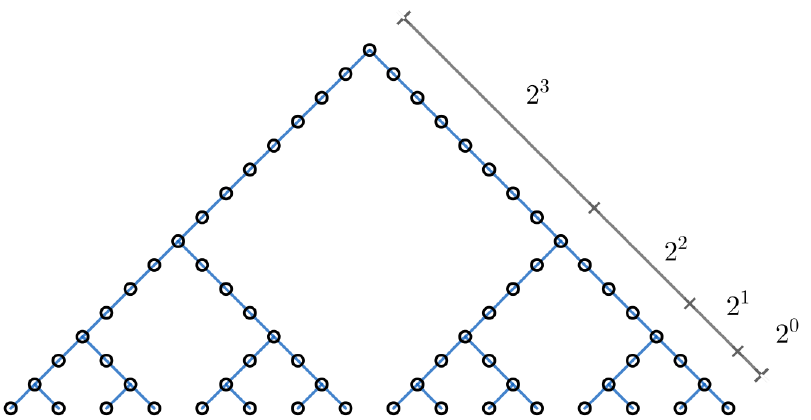}
    \caption{The   tree  $W_{2,4}$}\label{figure-W23}
\end{figure}

The tree $W_{\vt,n}$
is constructed by  replacing each edge of $B_n$  on the level $k$ by a path of length 
$\vt^{n-k}$, 
see Figure~\ref{figure-W23}. 
Thus the  
tree $W_{\vt,n}$ may be seen, as described in \cite{BMSZ21}, as a ``binary tree with variable size legs''.
Alternatively, one may think of  $W_{\vt,n}$ as a subgraph of the binary tree $B_m$, where $m=\sum_{k=1}^{n}\vt^{n-k}=(\vt^n-1)/(\vt-1)$ and multiple branches of $B_m$ are removed. Note that  $W_{\vt,n}$ is an unweighted tree
and  that $W_{2,n}$ is the same graph as the one described in the example in \cite[Section~7]{BMSZ21}.

Vertices of $W_{\vt,n}$  that came from the original tree $B_n$ will be called {\it nodes} or {\it branching points} (note that all nodes, except for the root and the leaves of the tree, have degree equal to 3, and all other vertices of 
$W_{\vt,n}$ have degree at most 2).  We  say that a node of $W_{\vt,n}$ is on level $k$, $0\le k\le n$, if it was a vertex of level $k$ in $B_n$.

For $n\in\bbN$ and $\vt\in\bbN$ with $\vt\ge 2$, the weighted tree $\WB_{\vt,n}$ is obtained by taking the complete binary tree $B_n$ and assigning to each edge of $B_n$ at level $k$ a weight equal to $\vt^{n-k}$.  Note that the sets of edges and vertices of  $\WB_{\vt,n}$ coincide with those of $B_n$ -- the only difference are the weights assigned to all nonterminal edges.

To simplify notation, when the values of $n$ and $\vt$ are clear from the context,  we will denote the metrics on $W_{\vt,n}$, $\WB_{\vt,n}$,  $\Lam(W_{\vt,n})$, and on   $\Lam(\WB_{\vt,n})$ by $d_W$, $d_{\WB}$,  $d_L$, and $d_{\LWB}$, respectively.

\section{Proof of Theorem~\ref{main1}}

We start from a simple fact that is surely  known. Since we could not find 
 it in the literature,
 for the convenience of the reader, we inlude its short proof.

\begin{lem} \label{KmPath}
  For any integer $m\ge 2$,   distortion of any embedding of $K_m$ into a path is at least $m-1$.
\end{lem}
\begin{proof}
    If $m$ points are mapped into a path, there are two which are furthest apart, and $m-2$ points between them. Denote by $d$ the distance between the first and last points. The remaining points partition this distance into $m-1$ subpaths. Thus the minimum distance between two of these points is at most $\frac{d}{m-1}$. Hence, the distortion of an embedding of $K_m$ is at least 
    ${d}/({\frac{d}{m-1}})=m-1$.
\end{proof}

\begin{theorem} \label{K4Embed}
  For any $n\in\bbN$  and $\vt\in\bbN$ with $\vt\ge 2$, the  distortion of any embedding of $K_4$ into $W_{\vt,n}$ or into  $\WB_{\vt,n}$ is at least $\frac{3}{2}$.
\end{theorem}
\begin{proof} Let $f:K_4\to W_{\vt,n}$ be a map onto the set $\{v_i\}_{i=1}^4\subset W_{\vt,n}$.
    By Lemma~\ref{KmPath}, if  3 or more points of this set are on a path, the distortion of $f$ is at least 2. So, assume that there is never 3 points on a path. Thus there must exist a branching node $s_1$ between $v_1$ and $v_2$ such that $s_1$ is on the geodesic paths  $p(v_3,v_1)$ and   $p(v_3,v_2)$. Further, there must exist a  a branching node $s_2$ on one of $p(s_1,v_i)$, $i=1,2,$ or $3$, such that $s_2$ is on $p(v_4,v_i)$ for each 
$i=1,2,3$. \Buo we suppose that $s_2$ is on $p(s_1,v_2)$ and
we denote, as marked on  Figure~\ref{fig:K4},
    $a_1=d_W(s_1,v_3)$, $a_2=d_W(s_1,v_1)$, $c_1=d_W(s_2,v_4)$, $c_2=d_W(s_2,v_2)$,  $a=\max\{a_1,a_2\}$,  $c=\max\{c_1,c_2\}$, and $b=d_W(s_1,s_2)$.
    \begin{figure}
        \centering
        \includegraphics[width=0.4\linewidth]{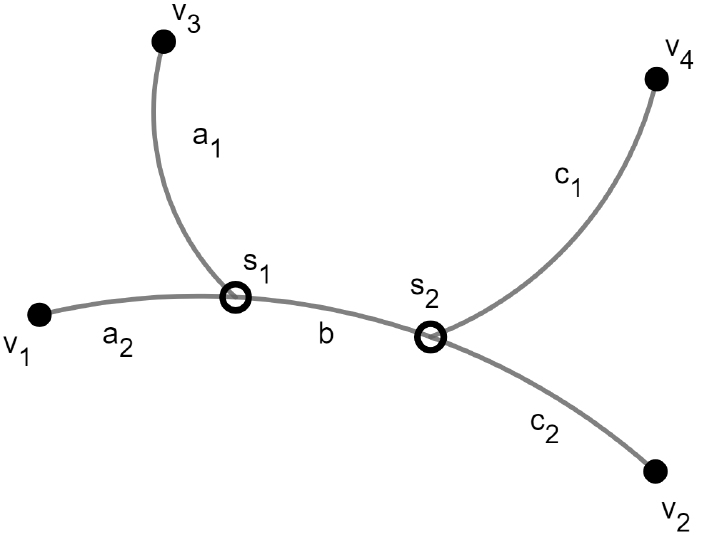}
        \caption{A possible configuration of an embedding of  $K_4$ into $W_{2,n}$, where $s_1$ and $s_2$ are branching points.}
        \label{fig:K4}
    \end{figure}

   \Buo we suppose that the branching level $k$ of $s_1$ is  greater than or equal to the branching level of $s_2$.
    Then $b\ge \vt^{n-k}$ and both vertices $v_1$ and $v_3$ have to be below vertex $s_1$, and  hence we get
   \begin{equation}\label{a}
   a=\max\{a_1,a_2\}\le \sum_{j=1}^{n-k} \vt^{n-(k+j)}=\frac{\vt^{n-k}-1}{\vt-1}<b.
    \end{equation}
This implies that 
 \begin{equation*}
  a_1+a_2\le \min\{a_1,a_2\} + \frac{\vt^{n-k}-1}{\vt-1}<\min\{a_1,a_2\} +b,
    \end{equation*}
and therefore the shortest distance between vertices  in the image of $K_4$ is at most 
\[\min\{ a_1+a_2, c_1+c_2\}\le      2\min\{a,c\}.\]
Since the longest distance  between vertices  in the image of $K_4$ is at least $a+b+c$ and, by \eqref{a}, 
$b\ge \min\{a,c\}$, we obtain 
    \[dist(f) \ge \frac{a+b+c}{2\min\{a,c\}}\geq 
\frac{3\min\{a,c\}}{2\min\{a,c\}}=\frac{3}{2},\]
which ends the proof for the case of $W_{\vt,n}$.

The case of the weighted trees  $\WB_{\vt,n}$ is very similar. The only difference is that in   $\WB_{\vt,n}$ every vertex is a branching node, but, due to the assigned  weights to edges at level $k$, the distance between vertices $s_1$ and $s_2$ is still at least $\vt^{n-k}$. Thus the rest of the proof is identical as in case of $W_{\vt,n}$.
\end{proof}

As a consequence of Theorem~\ref{K4Embed}, by a standard Ramsey theory argument, we obtain
a proof of Theorem~\ref{main1}.

\begin{proof}
    Suppose for contradiction that  there exists $C\in\mathbb{R}$ such that for all $k\in\mathbb{N}$ there exist $n\in\bbN$ and a 
    bi-Lipschitz embedding $\phi$ of $K_k$ into $W_{\vt,n}$ (respectively,  $\WB_{\vt,n}$) with distortion less than $C$.

    Let $m\in\bbN$ be the smallest integer such that $(3/2)^m>C$.
By the Ramsey's Theorem (see e.g. \cite[Theorem~3.6]{Rob2021})  there exists the Ramsey number $R(m,4)\in\bbN$ such that for every $k\ge R(m,4)$, every $m$-coloring of edges of $K_k$ contains a monochromatic copy of $K_4$.

Let $k> R(m,4)$, and  $n\in\bbN$ be such that there exists  an 
     embedding $\phi$ of $K_k$ into $W_{2,n}$ with distortion at most $C$.
Let $l,L>0$ be optimal numbers such that for all $u,v\in K_k$ we have
\[l\le d(\phi(u),\phi(v))\le L.\]
Then $L< Cl< (3/2)^m l.$ We color an edge $(u,v)$ of $K_k$ with color $i\in\{1,\dots,m\}$ iff
\[\Big(\frac32\Big)^{i-1}l\leq d(\phi(u),\phi(v))<\Big(\frac32\Big)^il.\]

Since $k>R(m,4)$, there exist a monochromatic copy  of $K_4\subset K_k$ with respect to this $m$-coloring. Since $\phi$ restricted to any monochromatic subset has distortion less than 
$3/2$, this contradicts Theorem~\ref{K4Embed}, and ends the proof.
\end{proof}

\section{Proof of Theorem~\ref{main2}}
\begin{figure} 
    \centering
    \includegraphics[width=0.5\linewidth]{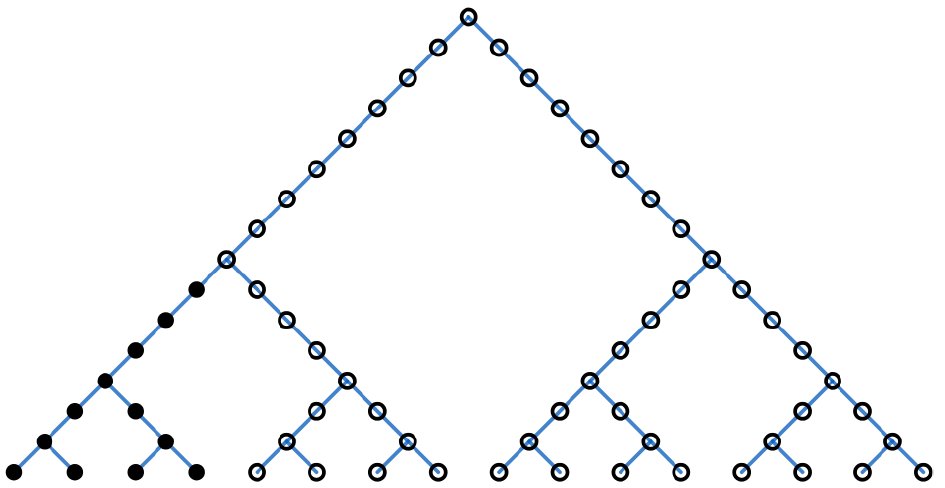}
    \caption{The set $F^4_{1,1}\subset W_{2,4}$}\label{fig-F-k}
\end{figure}

Let $n\in\bbN$ with $n\ge 3$.
For any $k\in \{0,\dots, n-1\}$, there are $2^k$ nodes of level $k$ in $W_{2,n}$, and there are 2 downward paths of length $2^{n-(k+1)}$ from each node. We index these paths from  left to right by 
an $i\in\{1,\dots, 2^{k+1}\}$. We define the set $F_{k,i}^{n}\subset W_{2,n}$ to be the set of all vertices that are in the $i$-th path and below, see 
Figure~\ref{fig-F-k}.
\begin{lem} \label{DkLemma}
For any $n\in \bbN$, $0\le k<n$, and $1\le i\le 2^{k+1}$, 
the cardinality of the set $F_{k,i}^{n}\subset W_{2,n}$ does not depend on $i$ and equals
\begin{equation}\lb{cardF}
 |F_{k,i}^{n}|:=w(n,k)=(n-k)2^{n-(k+1)}.
\end{equation}
\end{lem}
\begin{proof}
Indeed, the length of the path in $W_{2,n}$ from a node at level $k$ to a node at level $k+1$  below it is equal to $2^{n-(k+1)}$. The paths at the next level are of length $2^{n-(k+2)}$, but there are two of them.
For each  level, there are $2^{j-1}$ paths, each of length $2^{n-(k+j)}$. Thus 
\[|F_{k,i}^{n}|=\sum_{j=1}^{n-k} 2^{j-1}2^{n-(k+j)}=(n-k)2^{n-(k+1)}.\]
\end{proof}

 \begin{figure}
        \centering
        \includegraphics[width=0.5\linewidth]{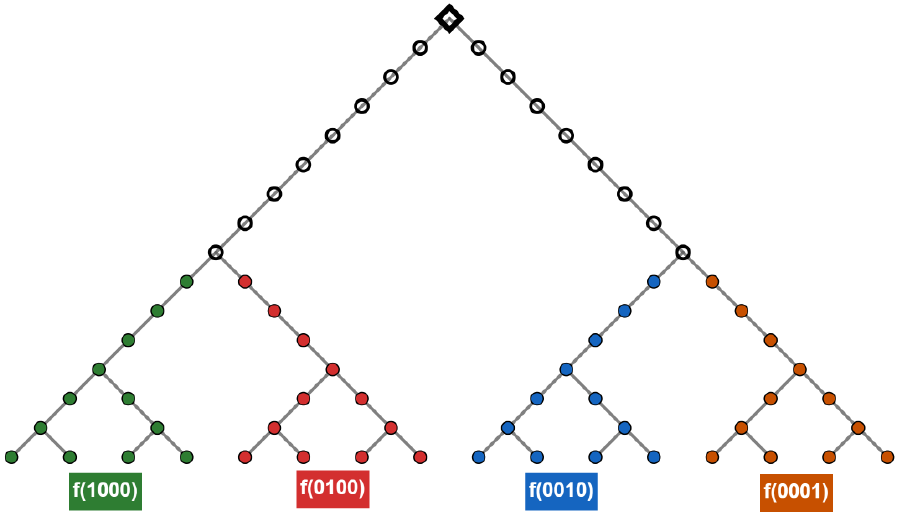}
        \caption {Images of the four basis elements of $H_4$  in $\Lam(W_{2,4})$ - in each of them the lampligher is at the root of $W_{2,4}$}\label{fig:H4}
    \end{figure}

 Theorem~\ref{main2} will follow immediately from our next result.

\begin{theorem} \label{H2kEmbed}
Let $k\in\mathbb{N}$ and $\e>0$. Then $H_{2^k}$  bi-Lipschitzly embeds into $\Lam(W_{2,n})$  with distortion at most $1+\e$, whenever $n>(2^{k+1}/3\e)+k-1$.
\end{theorem}

\begin{proof}[Proof of Theorem~\ref{H2kEmbed}]
Let $n>(2^{k+1}/3\e)+k-2$. We define the map  $f: H_{2^k}\to \Lam(W_{2,n})$ by setting for each $x\in H_{2^k}$
\begin{equation}\lb{embed}
f(x):=\big( \bigcup_{\{i : x_i=1\}} F_{k-1,i}^n , v_0\big),
\end{equation}
cf. Figure~\ref{fig:H4}.

Given any $x,y\in H_2^k$, let $I=\{ i\in [2^k] : x_i\ne y_i\}$. Then $|I|=d_H(x,y)$, and
\[d_L(f(x), f(y))=d_L\Big(\big( \bigcup_{i\in I} F_{k-1,i}^n , v_0\big), (\emptyset, v_0)\Big).\]

Note that 
a move in $\Lam(W_{2,n})$ from the element $(\emptyset, v_0)$ to
 $(\bigcup_{i\in I} F_{k-1,i}^n, v_0)$ requires that we travel from the root  to the top node of each of the sets $F_{k-1,i}^n$, $i\in I$, and back. 
 The distance from the root $v_0$ to any node of level $k-1$ is equal to 
\[\varrho_{k-1}=\sum_{j=1}^{k-1} 2^{n-j} = 2^{n}-2^{n-k+1}.
\]
The paths from the root to certain nodes may have significant overlaps depending on the relative positions of the nodes, but in any case, we have an upper estimate of the total travel distance to all nodes for $i\in I$ and back by 
$|I|\cdot 2\varrho_{k-1}$.

Also that for each $i\in I$, we  need to travel down  from the top node of   $F_{k-1,i}^n$,
$i\in I$,  switch on all its lamps, and return back to the top node.. As we computed in Lemma~\ref{DkLemma}, such travel in each set $F_{k-1,i}^n$ requires $3w(n,k-1)$ steps, and gives the total of $|I|\cdot 3w(n,k-1)$ steps to cover 
 all sets 
$F_{k-1,i}^n$, $i\in I$.

 Since $|I|=d_H(x,y)$,  we obtain 
\begin{equation*}
3w(n,k-1)d_H(x,y)\le
d_L(f(x), f(y))\le \big(3w(n,k-1)+2\varrho_{k-1}\big)d_H(x,y).
\end{equation*}

Therefore, by \eqref{cardF}, we estimate the distortion of the map $f$ as follows
\begin{equation*}
\begin{split}
dist(f)&\le \frac{3w(n,k-1)+2\varrho_{k-1}}{3w(n,k-1)}=1+\frac{2( 2^{n}-2^{n-k+1})}{3(n-k+1)2^{n-k}}\\
&
<1+\frac{2^{k+1}}{3(n-k+1)}<1+\e,
\end{split}
\end{equation*}
which ends the proof.
\end{proof}

\begin{rem}
It is easy to check that Theorem~\ref{main2} is valid also for the weighted trees $\WB_{2,n}$, with only minor changes in the proof. 
\end{rem}

\section{Proof of Theorem~\ref{nodiamonds}}

\begin{figure}
        \centering
        \includegraphics[width=0.4\linewidth]{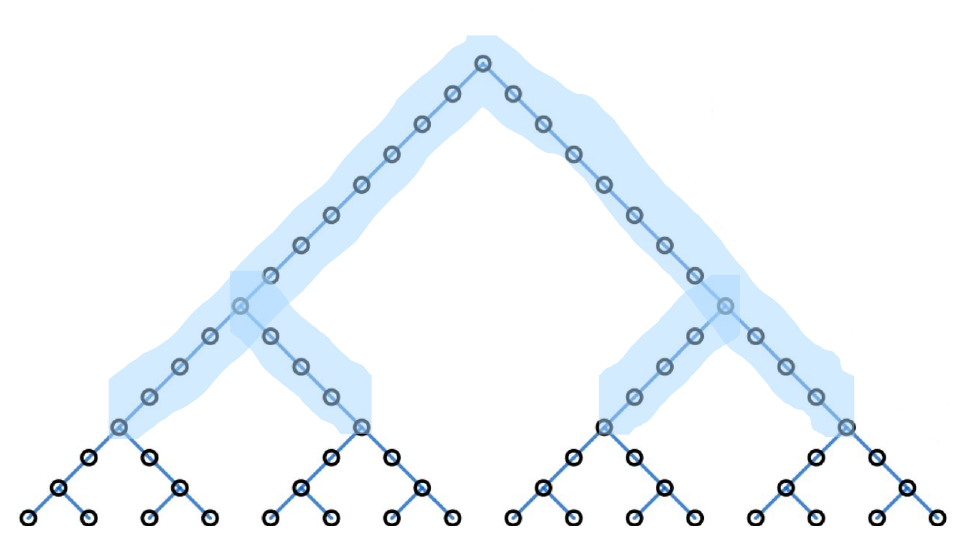}
        \caption {The set $F_n\subset W_{\vt,n}$  - the fork of depth 2}
        \label{fig:Fn}
    \end{figure}

We denote by $F_n$ the top ``fork of depth 2'' contained in $W_{\vt,n}$, see Figure~\ref{fig:Fn}.

\begin{lem} \lb{diam-lem}
For all $\vt, n\in\bbN$, with $\vt\ge 3$,
\begin{equation}\lb{diam}
\diam(\Lam(W_{\vt,n})\le  6\vt^n.
\end{equation}
\end{lem}

\begin{proof}
Since $W_{\vt,n}$ is an unweighted tree, by \cite[Theorem~3.1]{BMSZ21}, 
 we have $\diam(\Lam(W_{\vt,n})\le 3 |W_{\vt,n}|$. 
We estimate the cardinality of $W_{\vt,n}$ similarly as in the proof of \eqref{cardF}.
\begin{equation*}
\begin{split}
|W_{\vt,n}|&\le \sum_{j=1}^n 2^{j} \vt^{n-j}
=\vt^n \sum_{j=1}^n \Big(\frac2{\vt}\Big)^j\le \vt^n\cdot\frac{2}{\vt-2}\le 2\vt^n,
\end{split}
\end{equation*}
since $\vt\ge 3$, which ends the proof of \eqref{diam}.
\end{proof}

Stein and Taback \cite{ST13} showed that $\Lam(\bbZ)$ embeds bi-Lipschitzly with distortion at most 4 into a direct product of two infinite binary trees  $B_\infty$, cf. \cite[Proposition~2.2]{BMSZ21} for a different proof. 

Given two pointed graphs $G_i=(V_i, E_i,v_i)$ for $i=1,2$, we define the {\it vertex-coalescence $G_1\ast G_2$} of $G_1$ and $G_2$ by taking  $V(G_1\ast G_2)$ to be the disjoint union of $V_1$ and $V_2$ with  vertices $v_1$ and $v_2$ identified with each other to form one vertex denoted by $v_{0}$, and $E(G_1\ast G_2)$ is a disjoint union of $E_1$ and $E_2$. The graph metric on $G_1\ast G_2$ is given by
\begin{equation*}
d_{G_1\ast G_2}(u,v)=\begin{cases} d_{G_i}(u,v) &\text{ if } u,v\in G_i,\\
d_{G_1}(u,v_0) + d_{G_2}(v_0,v)  &\text{ if } u\in G_1 \text{ and } v\in G_2.
\end{cases}
\end{equation*}

For any $n\in \bbN$, the fork graph $F_n\subset W_{\vt,n}$ can be obtained by  coalescing three path graphs (at two different vertices),  one of length $2\cdot  \vt^{n-1}$ and two of length $2\cdot \vt^{n-2}$. It is well-known that for any finite $k\in\bbN$, the graph obtained by coalescing  at a common vertex $k$ copies of paths of arbitrary finite lengths embeds into the tree $B_\infty$ with distortion at most 2 (cf. \cite{Mat99}). Hence
 \cite[Theorem~4.2]{BMSZ21} shows that for any $n\in\bbN$, the lamplighter graph
  $\Lam(F_n)$ bi-Lipschitzly embeds with uniform distortion into a direct product of ten binary trees $B_\infty$.
Eskenazis \cite[Chapter~2]{E19} showed  that the sequence of diamonds
 $(\mD_k)_{k\in\bbN}$ does not equi-bi-Lipschitzly embed into any finite direct product of infinite binary trees $B_\infty$. The combination of these results gives 

\begin{cor}\lb{AE}
The sequence of diamonds $(\mD_k)_{k\in\bbN}$ does not equi-bi-Lipschitzly embed into 
$\Lam(\bbZ)$ or into the sequence $(\Lam(F_n))_{n\in\bbN}$.
\end{cor}

We are now ready for the proof of Theorem~\ref{nodiamonds}.

\begin{proof}
Fix $\vt\ge 3$.
Suppose, for contradiction, that there exists $C\in\bbN$, such that for all $k\in\bbN$, there exists $n\in\bbN$ with
\begin{equation*}
c_{\Lam(W_{\vt,n})}(\mD_k)\le 2^C.
\end{equation*}

By Corollary~\ref{AE}, there exists $k_C\in\bbN$ such that for any $n\in \bbN$,
\begin{equation}\lb{k0}
c_{\Lam(\bbZ)}(\mD_{k_C})> 2^C \text{ and } c_{\Lam(F_n)}(\mD_{k_C})> 2^C.
\end{equation}

Fix $k> k_C+2C +\vt+9$,  and let $n$ be the smallest natural number 
such that there exist a map $\vf:D_k\lra \Lam(W_{\vt,n})$ and $\la>0$ so that  for all $u,v\in \mD_k$,
\begin{equation}\lb{phik}
\la d_\mD(u,v)\le d_{L}(\vf(u), \vf(v)) \le 2^C \la d_\mD(u,v).
\end{equation}

Thus, by \eqref{diam} we have
\begin{equation*}
\la 2^k\le d_{L}(\vf(s), \vf(t)) \le 6\cdot \vt^n=6\vt^2\cdot \vt^{n-2}.
\end{equation*}
It is easy to check that if $\vt\ge 3$, then $6\vt^2\le 2^{\vt+3}$. Therefore
\begin{equation}\lb{la-est}
\la \le 2^{\vt+3-k} \vt^{n-2}.
\end{equation}
 \begin{figure}
        \centering
        \includegraphics[width=0.6\linewidth]{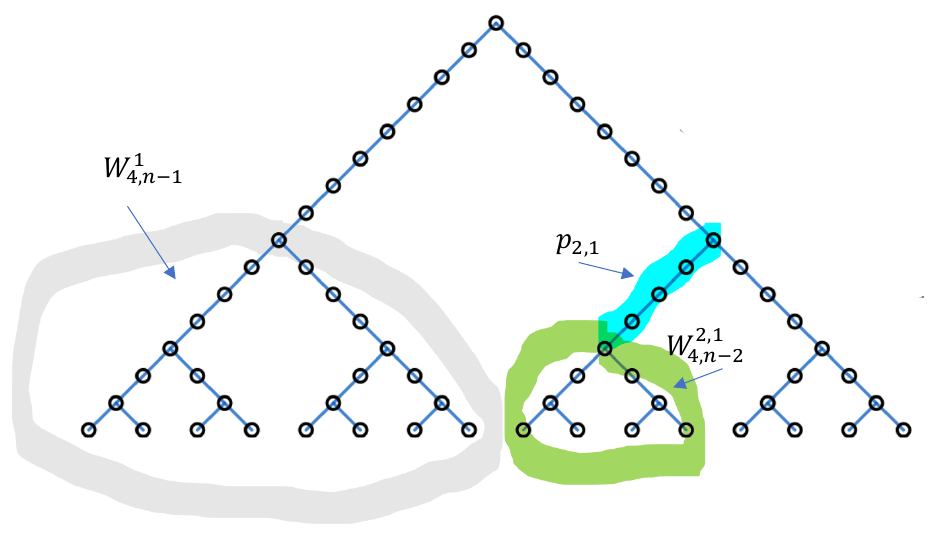}
        \caption {The subsets $ W_{\vt,n-1}^1,  W_{\vt,n-2}^{2,1}$, and $p_{2,1}$ of $W_{\vt,n}$.}
        \label{fig:Wn-1etc}
    \end{figure}

By \eqref{k0} and the minimality of $n$,
for any $\nu\in\{1,2\}$,
\[\vf(\mD_k)\not\subseteq \Lam(F_n)  \text{ and }  \vf(\mD_k)\not\subseteq \Lam(W^\nu_{\vt,n-1}).\]

Thus   there exist $u,v\in \mD_k$, with $\vf(u)=(a,A)$, $\vf(v)=(b,B)$, such that there exist $x, y \in \{a,b\}\cup (A\symdif B)$ and $\nu,\mu\in\{1,2\}$ with 
$x\in W^{\nu,\mu}_{\vt,n-2}$ and
 $y\not\in W^\nu_{\vt,n-1}$. Hence  $d_W(x,y)>\vt^{n-2}$ and
any walk in  from $x$ to $y$ $W_{\vt,n}$ has to contain the path
 $p=p_{\nu,\mu}$ from  
$t\DEF t_{\nu,\mu}$ to $\bar{t}\DEF t_\nu$ to . Note that the length of $p$ is equal to $\vt^{n-2}$.

 Define $z=z_{\nu,\mu},\bar{z}= \bar{z}_{\nu,\mu}$ to be the  points on $p$ with
\begin{equation}\lb{middlex}
d_W(z,\bar{z})= \frac{\vt^{n-2}}{2^{C+2}}\ \text{ and }\ 
d_W(t,z)=d_W(\bar{z},\bar{t})=\frac12\Big(1-\frac1{2^{C+2}}\Big)\vt^{n-2}\DEF \bar{d}.
\end{equation}
Note that, since $C\ge 1$, we have
\begin{equation}\lb{dbarest}
d_W(z,\bar{z})\le \frac{1}{2^C} \bar{d}.
\end{equation}

Let $u=u_0, u_1,\dots,u_m=v$ be a path in $\mD_k$ from $u$ to $v$, that is, $m=d_\mD(u,v)$ and for each 
$i\le m$, $(u_{i-1},u_i)$ is an edge in $\mD_k$.  Since $k-(C+\vt+3)\ge C+4$, 
by \eqref{la-est}, we have  for all $i$, 
\begin{equation}\lb{edge-est}
d_L(\vf(u_{i-1}), \vf(u_i)) \le 2^C\la \le 2^{C+\vt+3-k} \vt^{n-2}\le \frac14 d_W(z,\bar{z}).
\end{equation}
For each $i\le m$, let $\vf(u_i)=(A_i,a_i)$ and let $p_i$ denote a geodesic path from $\vf(u_{i-1})$ to 
$\vf(u_i)$ in  $ \Lam(W_{3,n})$.  The length of each $p_i$ is estimated from above in \eqref{edge-est}.
Since $w=\bigcup_{i=1}^m p_i$ is a 
 walk  in $ \Lam(W_{\vt,n})$ from $\vf(u)=(A,a)$ to $\vf(v)=(B,b)$, $w$ has to contain one or more elements
of each  of the forms  $(I,z)$ and  $(J,\bar{z})$, where $I,J\subseteq W_{\vt,n}$. 
Let $\al\in\{1,\dots, m\}$ be the largest index such that $(I,z)\in p_\al$ for some set $I$. 
Then $d_W(a_\al,z)$ is less than or equal to the length of the path $p_\al$, so,  
by \eqref{middlex},  \eqref{edge-est}, and \eqref{dbarest}, we obtain
\begin{equation}\lb{dista_al}
\begin{split}
\min\{d_W(a_\al,t),d_W(a_\al,\bar{t})\}&\ge \min\{d_W(t,z),d_W(\bar{t},z)\}-d_W(a_\al,z)\\
&\ge \bar{d} -  \frac14 d_W(z,\bar{z}) \ge  \bar{d} -  \frac1{2^{C+2}} \bar{d}\\
&>\frac34 \bar{d}.
\end{split}
\end{equation}

 Let 
$\be\in\{\al,\dots,m\}$ be the smallest index larger than or equal to $\al$ so that
 $(J,\bar{z})\in p_{\be}$ for some set $J$. 
By \eqref{edge-est}, we obtain
\begin{equation*}
\begin{split}
d_L(\vf(u_{\al}), \vf(u_\be)) &\ge d_L((I,z), (J,\bar{z}))-d_L((I,z),\vf(u_{\al}))-
d_L((J,\bar{z}), \vf(u_\be))\\
&\ge d_W(z,\bar{z})-length(p_\al)-length(p_\be)\\
&\ge d_W(z,\bar{z})-\frac{2}{4} d_W(z,\bar{z})\\
&= \frac{1}{2} d_W(z,\bar{z})
= \frac{1}{2^{C+3}}\vt^{n-2}.
\end{split}
\end{equation*}

Let $\g\in \{\al,\dots,\be\}$ be the smallest index such that 
\begin{equation}\lb{low-est}
d_L(\vf(u_{\al}), \vf(u_\g)) 
\ge \frac{1}{2^{C+3}}\vt^{n-2}=\frac{1}{2} d_W(z,\bar{z}).
\end{equation}

 Then, 
by \eqref{edge-est} and \eqref{dbarest},
\begin{equation}\lb{up-est}
\begin{split}
d_L(\vf(u_{\al}), \vf(u_\g)) 
&\le d_L(\vf(u_{\al}), \vf(u_{\g-1}))+d_L(\vf(u_{\g-1}), \vf(u_\g))\\
&\le \frac{1}{2} d_W(z,\bar{z}) +\frac{1}{4} d_W(z,\bar{z})\\
&= \frac{3}{4} d_W(z,\bar{z}) \le  \frac{1}{2^{C}} \cdot \frac{3}{4}\bar{d}.
\end{split}
\end{equation}

 By \eqref{phik}, \eqref{low-est},  and \eqref{la-est},   we get
\begin{equation*}
\begin{split}
d_\mD(u_\al,u_\g)&\ge \frac{1}{2^C\la}d_L(\vf(u_{\al}), \vf(u_\g)) \\
&\ge \frac{\vt^{n-2}}{2^{2C+3}2^{\vt+3-k}\vt^{n-2}}=2^{k-(2C+\vt+6)}.
\end{split}
\end{equation*}

Thus, since $k\ge 2C+\vt+9$, by Lemma~\ref{subdiamond}, 
 there exists a subdiamond  $\bar{\mD}=\mD_{k-(2C+\vt+9)}\subset \mD_k$ 
that lies between $u_\al$ and $u_\g$, that is, such that, for all $w\in\bar{\mD}$, 
\begin{equation}
\lb{up-est2}
d_\mD(u_\al,w)+d_\mD(w,u_\g)=d_\mD(u_\al,u_\g).
\end{equation}

Thus, by \eqref{up-est2}, \eqref{phik},  and \eqref{up-est}, for all $w\in\bar{\mD}$, 
\begin{equation}\lb{final-est}
\begin{split}
d_L(\vf(u_\al),\vf(w))+d_L(\vf(w),\vf(u_\g))
&\le 2^C\la(d_\mD(u_\al,w)+d_\mD(w,u_\g))\\
&= 2^C\la d_\mD(u_\al,u_\g)\\
&\le 2^C d_L(\vf(u_\al),\vf(u_\g))\\
&\le \frac{3}{4} \bar{d}.
\end{split}
\end{equation}

Recall that we denoted $\vf(u_\al)=(A_\al, a_\al)$ and  $\vf(u_\g)=(A_\g, a_\g)$. For each $w\in\bar{\mD}$ we denote $\vf(w)=(A_w, a_w)$. Then, by \eqref{tsp} and \eqref{final-est}, for $j\in\{\g,w\}$ we have
\begin{equation*}
\begin{split}
\max\Big\{ d_W(a_\al,a_j), \max_{s\in A_\al\symdif A_j}d_W(a_\al,s)\Big\}\le d_L((A_\al, a_\al),(A_j, a_j))
\le \frac{3}{4} \bar{d}.
\end{split}
\end{equation*}

Hence, by \eqref{dista_al}, the vertices $a_\g, a_w$ belong to the path $p$ connecting $t$ and $\bar{t}$,
and $(A_\al\symdif A_\g)\cup (A_\al\symdif A_w)\subset p$. Since 
$A_\g\symdif A_w \subset (A_\al\symdif A_\g)\cup (A_\al\symdif A_w)$, we also have 
$A_\g\symdif A_w \subset p$.
Thus $\vf_p:\bar{\mD}\lra  \Lam(W_{\vt,n})$ defined by $\vf_p(w)=(A_w\cap p,a_w)$ is a bi-Lipschitz embedding of 
$\bar{\mD}$ into $\Lam(p)$ with distortion at most $2^C$.
Since   $\bar{\mD}=\mD_{k-(2C+\vt+9)}$ and $k-(2C+\vt+9)\ge k_C$, this  contradicts \eqref{k0}, and thus ends the proof of Theorem~\ref{nodiamonds}.
\end{proof}


\begin{small}

\end{small}

\textsc{Department of Mathematics, Miami University, Oxford, OH
45056, USA} \par
  \textit{E-mail address}: \texttt{charlottemelby78@gmail.com} \par

\textsc{Department of Mathematics, Miami University, Oxford, OH
45056, USA} \par
  \textit{E-mail address}: \texttt{randrib@miamioh.edu} \par

\end{document}